\documentclass[12pt]{amsart}
\usepackage[T2A]{fontenc}
\usepackage[cp1251]{inputenc}
\usepackage{amsmath,amsfonts,textcomp}
\usepackage{amssymb}
\usepackage{amscd}
\usepackage[dvips]{graphicx}

\renewcommand{\Re}{\mathop{\rm Re}\nolimits}
\renewcommand{\Im}{\mathop{\rm Im}\nolimits}
\newcommand{\im}{\ensuremath{\mbox{Im}\,}}
\newcommand{\re}{\ensuremath{\mbox{Re}\,}}

\newcommand{\lemma}[1]{\bf Lemma #1.\it\,}
\newcommand{\proposition}[1]{\bf Proposition #1.\it\,}
\newcommand{\cor}[1]{\bf Corollary #1.\it\,}
\newcommand{\theor}[1]{\bf Theorem #1.\it\,}
\newcommand{\doc}{\bf Proof.\rm\,}
\newcommand{\remark}[1]{\bf Remark #1.\rm\,}
\newcommand{\CC}[1]{\mathbb{C}^{#1}}
\newcommand{\RR}[1]{\mathbb{R}^{#1}}
\newcommand{\dw}{\frac{\partial}{\partial w_1}}
\newcommand{\dww}{\frac{\partial}{\partial w_2}}
\newcommand{\dwww}{\frac{\partial}{\partial w_3}}
\newcommand{\dz}{\frac{\partial}{\partial z_1}}
\newcommand{\dzz}{\frac{\partial}{\partial z_2}}
\newcommand{\dzzz}{\frac{\partial}{\partial z_3}}
\newcommand{\z}{z_1}
\newcommand{\zz}{z_2}
\newcommand{\zzz}{z_3}
\newcommand{\w}{w_1}
\newcommand{\ww}{w_2}
\newcommand{\www}{w_3}
\newcommand{\xo}{X_0}
\newcommand{\x}{X_1}
\newcommand{\xp}{X'_1}
\newcommand{\xx}{X_2}
\newcommand{\xxx}{X_3}
\newcommand{\lr}{\longrightarrow}
\newcommand{\rr}{\Rightarrow}
\newcommand{\dd}{\partial}
\newcommand{\la}{\lambda}
\newcommand{\MP}[1]{M^{+}_{#1}}
\newcommand{\MM}[1]{M^{-}_{#1}}

\title{Homogeneous hypersurfaces in $\CC{3}$, associated with a model
 CR-cubic}
\author {V.K.Beloshapka and I.G.Kossovskiy}

\begin{document}
\maketitle

\begin{abstract}
The model 4-dimensional CR-cubic in $\CC{3}$ has the following
"model" property: it is (essentially) the unique locally homogeneous 4-dimensional
CR-manifold in $\CC{3}$ with finite-dimensional infinitesimal automorphism
algebra $\mathfrak{g}$ and non-trivial isotropy subalgebra. We study and classify, up to
local biholomorphic equivalence, all $\mathfrak{g}$-homogeneous hypersurfaces in
$\CC{3}$ and also classify the corresponding local transitive actions of the model algebra
$\mathfrak{g}$ on hypersurfaces in $\CC{3}$.
\end{abstract}

\section{introduction}

The most interesting objects in CR-geometry  are CR-manifolds with
symmetries, i.e. CR-manifolds, admitting (local) actions of Lie
groups by holomorphic transformations. If such an action is (locally)
transitive, then the manifold is called \it  (locally) holomorphically
homogeneous \rm (or just homogeneous). Locally homogeneous manifolds are
"the same in all points", i.e. the germs  of a locally
homogeneous manifold  at any two points   are biholomorphically equivalent. Among all
homogeneous CR-manifolds one can single out so-called model
manifolds - algebraic CR-submanifolds in $\CC{N}$ with maximal-dimensional
automorphism groups. As it was demonstrated in
\cite{chern},\cite{tanaka},\cite{obzor}, the properties of model
manifolds determine in many aspects the properties of general
CR-manifolds. In this paper some interplay between a model
4-manifold in $\CC{3}$ and homogeneous hypersurfaces in $\CC{3}$
is studied.

To work with homogeneous CR-manifolds and their symmetry groups
and algebras we give a few definitions.

 Consider in the complex space $\CC{N}$ a germ $M_p$ of a generic
real-analytic CR-submanifold $M$ at a point $p\in M$ (we
suppose that all CR-submanifolds are real-analytic and generic if not
otherwise mentioned). We  consider the following objects:

1) $\it \mathfrak{aut} M_{p}$ - the Lie algebra of germs at the point $p$ of vector fields
 of the form
$$2\re\left(f_1(z)\dz+\dots+f_N(z)\frac{\partial}{\partial z_N}\right),$$ which are
tangent to $M$ at each point, and the functions $f_j(z)$ are
holomorphic
 in a neighborhood of $p$. We  call such vector fields \it
holomorphic vector fields on $M$ in a neighborhood of $p$. \rm Clearly these vector fields
are exactly that ones which generate local actions of Lie groups
on $M$ by transformations, holomorphic in a neighborhood of $p$
in $\CC{N}$. The Lie algebra $\it \mathfrak{aut} M_p$ is called \it the
infinitesimal automorphism algebra of $M$ at $p$. \rm If $\it \mathfrak{aut}
M_{p}$ is finite-dimensional, then all vector fields from this
algebra can be defined in the same neighborhood and there exists a connected simply-connected
Lie group, acting on $M$ locally by holomorphic
transformations in a neighborhood of the point $p$, such that its
tangent algebra is isomorphic to $\it \mathfrak{aut} M_{p}$ and the vector
fields from $\it \mathfrak{aut} M_{p}$ are the infinitesimal generators of
the action. We  denote this local group by $\it Aut M_{p}$ and
call it \it the local holomorphic automorphism group of $M$ at $p$.
\rm

2) $\it \mathfrak{aut}_{p}M_{p}$ - the Lie subalgebra in
$\it \mathfrak{aut} M_{p}$, which consists of germs of vector fields from
$\it \mathfrak{aut} M_{p}$, vanishing at $p$. This algebra is called \it
the stability subalgebra of $M$ at the point $p$ or the isotropy
subalgebra. \rm If $\it \mathfrak{aut} M_{p}$  is finite-dimensional, then
$\it \mathfrak{aut}_{p}M_{p}$ is naturally identified with the tangent
algebra of the \it stability group \rm $\it Aut_{p}M_{p}$,
 which consists of holomorphic automorphisms of the germ, fixing the point $p$.

\it A local action \rm of a finite-dimensional real Lie algebra
$\mathfrak{h}$ on the germ $\CC{N}_p$ of the complex space
$\CC{N}$ at a point $p$ is a homomorphism $\varphi:\,\mathfrak{h}\lr\mathfrak{aut}\,\CC{N}_p$. If $M$ is a
CR-submanifold in $\CC{N}$, passing through $p$, we say that \it
$\mathfrak{h}$ acts transitively on $M_p$ (or that $\mathfrak{h}$
acts locally transitively on $M$ at the point $p$), \rm if the
linear space, spanned by the values at the point $p$ of the vector
fields from $\varphi(\mathfrak{h})$, which are tangent to $M$,  coincides
with $T_p M$. The germ $M_p$ is called \it homogeneous \rm in this
case, and the manifold $M$ is called \it locally homogeneous at
$p$. \rm If $M$ is locally homogeneous at all points, then we call
it just \it locally homogeneous. \rm For more information about possible
equivalent definitions of homogeneous CR-manifolds we refer to
\cite{zaitsev}.

With any local action of a Lie algebra $\mathfrak{h}$ we can
associate a local action of a local Lie group $H$ with tangent
algebra $\mathfrak{h}$ and consider the orbits of this action.
These orbits are locally homogeneous CR-manifolds and their local
homogeneity is provided by the Lie algebra $\mathfrak{h}$. We call
this collection of orbits \it locally homogeneous manifolds, associated
with the Lie algebra $\mathfrak{h}$. \rm

Coming back to homogeneous CR-manifolds in $\CC{3}$, we  firstly
mention E.Cartan's classification theorem for homogeneous
hypersurfaces in $\CC{2}$ (see \cite{cartan}). Due to this
theorem, the following trichotomy holds for a locally
homogeneous hypersurface in $\CC{2}_{z,w}$:

(1) $\mbox{dim}\,\mathfrak{aut}M_{p}=\infty,$  which occurs if and only if
$M$ is locally biholomorphically equivalent to the hyperplane $\im
w=0$ (the Levi-flat case).

(2) $\mbox{dim}\, \mathfrak{aut}M_{p}=8,$ which occurs if and only if $M$ is
locally biholomorphically equivalent to the unit sphere
$|z|^2+|w|^2=1$.

(3) $\mbox{dim}\, \mathfrak{aut}M_{p}=3,$ which occurs if and only if $M$ is
locally biholomorphically equivalent to one of \it Cartan's
homogeneous surfaces \rm (see \cite{cartan} for details).

Note that due to a classical  result of H.Poincare \cite{poincare}, all
other hypersurfaces in $\CC{2}$ have infinitesimal automorphism
algebras of dimensions $\leq 2$.

Hence we have the following \it rigidity phenomenon \rm for germs
of Levi non-degenerate homogeneous hypersurfaces in $\CC{2}$: any
such germ is either a germ of the model surface (i.e. the sphere
in our case) and has maximal-dimensional infinitesimal
automorphism algebra, or it is \it holomorphically rigid, \rm i.e
its stability subalgebra is trivial.

The classification of homogeneous hypersurfaces in $\CC{3}$ is not
complete yet. In the case of Levi non-degenerate
hypersurfaces with high-dimensional isotropy subalgebras the
classification was obtained by A.Loboda (see
\cite{loboda1},\cite{loboda2}). In the Levi degenerate case the
full classification was obtained by G.Fels and W.Kaup. To describe
their results we  give the following definition: a
CR-submanifold $M$ in $\CC{N}$ is called \it holomorphically
degenerate, \rm if in a neighborhood  of any point $p\in M$ there exists a non-zero holomorphic vector field on $M$,
which belongs to the complex tangent space of $M$ at each point.
In this case it is not difficult to see that
$\mbox{dim}\, \mathfrak{aut}\rm M_{p}=\infty.$ Otherwise $M$ is called \it
holomorphically non-degenerate. \rm In particular, all Levi
non-degenerate hypersurfaces are holomorphically non-degenerate.
In case of Levi degenerate hypersurface in $\CC{3}$ this
non-degeneracy condition is equivalent to the \it 2-nondegeneracy,
\rm in a general point, which is some condition on the defining function of the
hypersurface (see \cite{bouendy} for details). For a
2-nondegenerate hypersurface the Levi form has rank 1 at each point and $\mbox{dim}\,\it \mathfrak{aut}M_{p}<\infty.$

Now we can formulate G.Fels and W.Kaup's classification theorem.
Due to this theorem, the following trichotomy holds for a locally homogeneous Levi
degenerate hypersurface in $\CC{3}$:

(1) $\mbox{dim}\, \mathfrak{aut}M_{p}=\infty,$  which occurs if and only
if $M$ is locally biholomorphically equivalent to a direct product
$M^3\times \CC{}$, where $M^3\subset\CC{2}$ is  one
of the homogeneous hypersurfaces in $\CC{2}$ from E.Cartan's
list specified above (holomorphically degenerate case).

(2) $\mbox{dim}\, \mathfrak{aut}M_{p}=10,$ which occurs if and only if $M$ is
locally biholomorphically equivalent to the tube over the future
light cone: $y_3^2=y_1^2+y_2^2,y_3>0,\,y_j=\im z_j,\,(z_1,z_2,z_3)\in\CC{3}$.

(3) $\mbox{dim}\, \mathfrak{aut}M_{p}=5,$ which occurs if and only if $M$ is
locally biholomorphically equivalent to the tube over an affinely
homogeneous hypersurface in $\RR{3}$ from some list, specified in
\cite{kaup}.

Hence,  in the same way as in E.Cartan's case, we have the
rigidity phenomenon  for germs of holomorphically non-degenerate locally
homogeneous hypersurfaces in $\CC{3}$: any such germ is either a
germ of the model surface (i.e. the tube over the future light
cone in that case) and has maximal-dimensional infinitesimal
automorphism algebra, or it is  holomorphically rigid, \rm i.e
its stability subalgebra is trivial.

We  study the class of  locally homogeneous hypersurfaces in $\CC{3}$
with the following property: the local homogeneity of these
surfaces is provided by one of the model algebras in $\CC{3}$ -
the unique 5-dimensional model algebra for the class of
4-dimensional holomorphically non-degenerate (or, equivalently, \it totally non-degenerate \cite{bes}\rm) CR-manifolds in
$\CC{3}$. This algebra is the infinitesimal automorphism algebra
$\mathfrak{g}$ of the model 4-dimensional CR-cubic $C$, given by
the following equations: $$\im w_2=|z|^2,\,\im w_3=2\re
(z^2\bar{z}),\,(z,w_2,w_3)\in \CC{3}$$
(this notation is related to a natural gradation of the coordinates in $\CC{3}$, see section 2).

Due to V.K.Beloshapka, V.V.Ezhov and G.Schmalz (see \cite{bes}),
the model properties of the cubic $C$ are given by the following trichotomy for a 4-dimensional locally homogeneous
 CR-manifold $M$ in $\CC{3}$:

(1) $\mbox{dim}\,\mathfrak{aut} M_p=\infty$, which occurs if and only if $M$
is locally biholomorphically equivalent to a direct product
$M^3\times \RR{1},\,M^3\subset\CC{2}_{z_1,z_2},\,\RR{1}\subset\CC{1}_{z_3}$ (the holomorphically degenerate
case).

(2)  $\mbox{dim}\,\mathfrak{aut} M_p=5$, which occurs if and only if $M$ is locally biholomorphically
equivalent to the cubic $C$.

(3) $\mbox{dim}\,\mathfrak{aut} M_p=4,\,\mbox{dim}\,\mathfrak{aut}_p M_p=0$ for all other manifolds (the rigidity phenomenon). \rm

It is proved also in \cite{bes} that the cubic $C$ is the most symmetric holomorphically non-degenerate 4-manifold in $\CC{3}$:
$\mbox{dim}\,\mathfrak{aut} M_0\leq \mbox{dim}\,\mathfrak{aut} C_0$, and the equality holds only for
 manifolds, locally biholomorphically equivalent to the cubic. The automorphism group $G$ and the infinitesimal automorphism
algebra $\mathfrak{g}$ of the cubic are described in the next section.

We  associate with the cubic $C$ some (locally) homogeneous
hypersurfaces in $\CC{3}$ in two different ways.

The first one is to consider the natural action of the
5-dimensional polynomial transformation group $G$ (or,
equivalently, of the model algebra $\mathfrak{g}$) in the ambient
space $\CC{3}$. Since the group is of dimension 5, we conclude
that the cubic is a singular 4-dimensional orbit of this action,
but general orbits are of dimension 5. This approach was realized
in \cite{orbits}. \it Note that according to the above trichotomy
for a homogeneous 4-manifold in $\CC{3}$ this machinery for the
construction of homogeneous hypersurfaces in $\CC{3}$, associated
with a homogeneous 4-manifold, is the only possible, i.e. the
obtained in \cite{orbits} class of hypersurfaces is (essentially)
the class of all locally homogeneous hypersurfaces, associated in
the natural sense with locally homogeneous 4-manifolds in $\CC{3}$
(in the category of non-degenerate manifolds). \rm In section 2 we
give a review of the results of \cite{orbits} and also give
another (tube) realization to the obtained in \cite{orbits}
foliation to orbits (case $A$ in the theorem below). It helps us
to recognize one special orbit as one of the hypersurfaces from
\cite{kaup} and also helps us to find an interesting realization
of the cubic $C$ as the tube over the twisted cubic in $\RR{3}$.
In section 3 we classify the obtained homogeneous hypersurfaces
and compute their automorphism groups. In particular, we prove an
analogue of the Poincare-Alexander theorem (see
\cite{poincare},\cite{alexander}) for the orbits under
consideration.

The second one is to consider all homogeneous hypersurfaces in
$\CC{3}$, associated with the abstract model Lie algebra
$\mathfrak{g}$. The approach in this case is analogue to
E.Cartan's approach in the classification problem for
hypersurfaces in $\CC{2}$. We find all possible realizations of the
abstract Lie algebra $\mathfrak{g}$ as an algebra of holomorphic
vector fields in $\CC{3}$, acting transitively on hypersurfaces,
and thus find  all possible orbits of the corresponding actions -
they  form the desired class of homogeneous hypersurfaces in
$\CC{3}$ (we call these hypersurfaces \it
$\mathfrak{g}$-homogeneous). \rm This approach  is realized in
section 4. In section 5 we  classify the obtained homogeneous
hypersurfaces up to local biholomorphic equivalence and compute
their infinitesimal automorphism algebras (and hence the
corresponding local automorphism groups). As a result we prove the
following classification theorem for $\mathfrak{g}$-homogeneous
hypersurfaces in $\CC{3}$:

\bf Main theorem. \it \noindent(1) The model algebra
$\mathfrak{g}$ has 4  types  of  local
transitive actions on hypersurfaces in $\CC{3}$ - actions of type
$A,A_1,A_0$ and $B$, described in section 5. Any two actions of different types are inequivalent. The corresponding
orbits look as follows.

\begin{gather}\notag  \mbox{TYPE}\,\, A:\,
 N_{\nu}^{\pm}=\left\{(y_3-3y_1y_2+2y_1^3)^2=\pm\nu(y_2-y_1^2)^3,\,\pm(y_2-y_1^2)>0\right\},\,\nu\geq 0,\\
\notag N^0=\left\{y_2= y_1^2, \,y_3\neq y_1^3  \right\}. \\
 \notag  \mbox{TYPE}\,\, A_1:\,\,\,
 S_\gamma=\left\{y_3=\gamma
y_1^3+\frac{y_2^2}{y_1},y_1\neq
0\right\},\,\gamma\in\RR{}.\\
\notag  \mbox{TYPE}\,\, A_0:\,\,\,  Q_\beta=\{y_3=\beta y_1^3+2y_2x_1,\,y_1\neq 0\},
\,\beta\in\RR{}.\\
\notag \mbox{TYPE}\,\, B:\,\,\, \Pi_\delta=\{y_2=\delta y_1,\,y_1y_2\neq
0\},\,\delta\in\RR{*}.\end{gather}

Here $z_k=x_k+iy_k$.

\noindent(2) Any $\mathfrak{g}$-homogeneous hypersurface in
$\CC{3}$ is locally biholomorphically equivalent to one of the
following pairwise non-equivalent homogeneous hypersurfaces in $\CC{3}$:

\noindent(a) Tube manifolds $N^\pm_{\nu}$ for $\nu>0$.

\noindent(b) Tube manifolds $S_{\pm 1}$ (the case of
$S_{\gamma},\gamma\in\RR{*})$.

\noindent(c) The tube over the future light cone
$y_3^2=y_1^2+y_2^2,\,y_3>0$ (the case of $S_0$).

\noindent(d) The indefinite quadric $y_3=|z_1|^2-|z_2|^2$ (the case of
$N^{\pm}_0$ and $Q_{\beta},\beta\in\RR{}$).

\noindent(e) The cylinder over the unit sphere in $\CC{2}$:
$|z_1|^2+|z_2|^2=1$ (the case of $M^0$).

\noindent(f) The real hyperplane $y_3=0$ (the case of
$\Pi_{\delta},\delta\in\RR{*}$).

\noindent In cases (e) and (f) the infinitesimal automorphism
algebras of the surfaces are infinite-dimensional; in cases (c)
and (d) these algebras are well-known simple Lie algebras (see
\cite{kaup1},\cite{obzor},\cite{chern} for the description of the
algebras and the  corresponding local automorphism groups); in
case (a) the infinitesimal automorphism algebras coincide with the
model algebra $\mathfrak{g}$, and all local automorphism of the
surfaces are global and belong to the group $G$; in case (b) the
infinitesimal automorphism algebras are  isomorphic to the model
algebra $\mathfrak{g}$  (more precisely, they coincide with the
algebra $A_1$), the corresponding local automorphism group is
described in section 5, hence in cases (a) and (b) the
hypersurfaces are holomorphically rigid.\rm

\remark{1.1} We note some interesting facts, which follow from
the above classification theorem.

(1) All $\mathfrak{g}$-homogeneous hypersurfaces are locally
biholomorphically equivalent to globally homogeneous
hypersurfaces.

(2) All $\mathfrak{g}$-homogeneous hypersurfaces turn out to be tube manifolds
over affinely homogeneous hypersurfaces in $\RR{3}$ in an appropriate local
coordinate system (for the indefinite quadric we get the tube realization by means of a quadratic variable change,
as well as for the unit sphere in the Poincare realization). Affinely homogeneous hypersurfaces in $\RR{3}$ were classified
in \cite{ezheas} and \cite{dcr}, but the corresponding tube manifolds in $\CC{3}$ were not studied from the point of view
of  holomorphic classification and automorphism groups.   Hence the present work can be considered as a step in this direction.

(3) The hypersurfaces $N^+_{\mu}$ for $\mu>0$, $N^-_{\nu}$ for
$\nu>0,\nu\neq 4$ and $S_{\pm 1}$ are Levi non-degenerate and
holomorphically rigid. Hence they give examples of pairwise
non-equivalent locally homogeneous hypersurfaces in $\CC{3}$,
which  are not covered by the classification theorems obtained in
\cite{loboda1},\cite{loboda2},\cite{kaup} (the exceptional orbit
$N^-_4$ is 2-nondegenerate and hence occurs in \cite{kaup}).

(4) In the same way as it results  in both E.Cartan's and
 G.Fels-W.Kaup's classification theorems, the following
rigidity phenomenon holds: each holomorphically non-degenerate
homogeneous hypersurface, generated in the specified sense by the
model algebra $\mathfrak{g}$, is either extremely-symmetric (a
quadric - the most symmetric Levi non-degenerate hypersurface, or
the tube over the future light cone - the most symmetric
2-nondegenerate hypersurface), or it is holomorphically rigid,
i.e. its isotropy subalgebra is trivial. Each of the obtained
infinitesimal automorphism algebras turns out to be isomorphic to
one of the model algebras in $\CC{3}$ (i.e. to the infinitesimal
automorphism algebra of a quadric, of the cubic or of the tube
over the future light cone). Thus the construction of homogeneous
hypersurfaces, used in \cite{orbits} and in the present paper,
gives an interesting connection among model algebras in $\CC{3}$.
It is also amazing that the obtained holomorphically degenerate
hypersurfaces are in a certain sense also extremely-symmetric: the
first one (the hyperplane) is the cylinder over the most symmetric
hypersurface in $\CC{2}$ - the hyperplane $\im w=0$, and the
second one is the cylinder over the most symmetric Levi
non-degenerate hypersurface in $\CC{2}$ - the unit sphere
$|z_1|^2+|z_2|^2=1$.

\remark{1.2} Note that the above classification theorem gives a
description of all possible hypersurface-type left-invariant
CR-structures on the group $G=Aut(C)$ (see section 2).

The authors would like to thank  W.Kaup for useful remarks, which helped to improve the
text of this paper.

\section{Action of the automorphism group of the cubic in the ambient space}

In this section we  describe the automorphism group $G$ and the
infinitesimal automorphism algebra $\mathfrak{g}$ of the cubic $C$,
and then give a review of the paper \cite{orbits}, where the
action of the group $G$  in the ambient
space $\CC{3}$ was studied and the corresponding orbits were
presented explicitly. Also we  present another (tube) realization
of the obtained foliation of $\CC{3}$, given by the group $G$,
which helps us to recognize one special orbit as a well-known
hypersurface in $\CC{3}$ and find a tube realization for the
model cubic $C$.

As it was mentioned in the introduction, the cubic $C$ is a
homogeneous 4-dimensional CR-manifold in $\CC{3}$, given by the
following equations: $$\im w_2=|z|^2,\,\im w_3=2\re
(z^2\overline{z}),\,(z,w_2,w_3)\in \CC{3}.$$ This notation is
associated with the following natural gradation of the coordinates
in $\CC{3}$:
\begin{gather}[z]=1,[w_2]=2,[w_3]=3.\end{gather} The polynomials
$\im w_2-|z|^2$ and $\im w_3-2\re (z^2\overline{z})$ are
homogeneous under this gradation and hence the cubic admits an
action of the following group of dilations: \begin{gather} z\lr\la
z,\,w_2\lr \la^2 w_2,\,w_3\lr \la^3
w_3,\,\la\in\RR{*}.\end{gather}

This group is the isotropy subgroup $G_0$ of the origin in the 5-dimensional group $G=Aut(C)$.
$G$ is a semidirect product of $G_0$ and the following
polynomial group $G_{-}$, providing the homogeneity of the cubic:
\begin{gather} \nonumber
z \mapsto  z+p,\\
w_2 \mapsto  w_2 + 2i\bar{p}  z + i|p|^2 + q_2,\\
\nonumber w_3 \mapsto  w_3 + 4(\Re p)  w_2 + 2 i( 2 |p|^2 +
\bar{p}^2)  z  + 2 i \bar{p} z^2 + 2 i \Re (p^2\bar{p})+q_3,\notag\\
\mbox{where}\,\,\, p\in\CC{},\,q_j\in\RR{}. \notag\end{gather}

The infinitesimal automorphism algebra $\mathfrak{g}$ of the
cubic, which can be naturally identified with the tangent algebra
of $G$, is a graded Lie algebra of kind
$$\mathfrak{g}=\mathfrak{g}_{-3}+\mathfrak{g}_{-2}+\mathfrak{g}_{-1}+\mathfrak{g}_{0},$$ where the gradation for
monomials is taken from (1), and the basic differential operators
are graded in the following way:
$$\left[\frac{\partial}{\partial
z}\right]=-1,\left[\dww\right]=-2,\left[\dwww\right]=-3.$$

The basic vector fields from $\mathfrak{g}$ look as follows (we
 skip the operator $2\re(\cdot)$):
\begin{gather*}
X_{3}=\dwww,\,X_2=\dww,\\
\xp=i\frac{\partial}{\partial
z}+2z\dww+2z^2\dwww,\\
\x=\frac{\partial}{\partial z}+2iz\dww+(4w_2+2iz^2)\dwww,\\
\xo=z\frac{\partial}{\partial z}+2w_2\dww+3w_3\dwww.\end{gather*}
Here $\mathfrak{g}_{-3}$ is spanned by $X_3,\,\,\mathfrak{g}_{-2}=<X_2>,\,\mathfrak{g}_{-1}=<X_1,\xp>,\,\mathfrak{g}_{0}=<X_0>$. Since $\mathfrak{g}_{0}$ is abelian, the algebra $\mathfrak{g}$ is
solvable. Also note that $$<X_3,X_2,\xp>$$ is an abelian ideal in
$\mathfrak{g}$.

The subalgebra $\mathfrak{g}_{0}$ is the isotropy subalgebra of
the origin and hence corresponds to the subgroup $G_0$, the
nilpotent ideal
$\mathfrak{g}_{-}=\mathfrak{g}_{-3}+\mathfrak{g}_{-2}+\mathfrak{g}_{-1}$
corresponds to the subgroup $G_{-}$ (this ideal coincides with the
unique irreducible 4-dimensional  nilpotent real Lie algebra).

The natural action of $G$ in the ambient space is given as a composition of actions (2) and (3).
Note, that the polynomial $P=\Im w_2-|z|^2$ is a relative
invariant of this natural action of weight 2 (more precisely, each
transformation from $G$ multiplies it by $\lambda^{2}$). Hence we
have 3 kinds of orbits: those lying in the domain $\Im
w_2-|z|^2>0$ (case 1 - orbits "over the ball"), those lying in the
domain $\Im w_2-|z|^2<0$ (case 2 - orbits "over the complement to
the ball"), and those lying over the quadric $\Im w_2-|z|^2=0$ (case
3 - orbits "over a sphere").

\it CASE 1. \rm In this case, as was shown in \cite{orbits}, for a
point $(a,b,c), \im b>|a|^2$ we get the following orbits:

\begin{gather*}
\Im w_3 = - 2 \Re z^2 \bar{z} + 4 \Re z \Im w_2 + |\mu| (\Im
w_2-|z|^2)^{\frac{3}{2}},\,\im w_2>|z|^2,\mu\in\RR{}.\end{gather*}

Any orbit  is an open smooth part of the real-analytic set
\begin{gather*}
\nonumber (\Im w_3 - 4 \Re z \Im w_2 + 2 |z|^2 \Re z)^2 =  \mu^2
(\Im w_2-|z|^2)^{3},
\end{gather*}
lying over $P>0$.

Any such orbit except the one with $\mu=0$ has two connected
components, corresponding to two $\mu$ with opposite signs. They
can be mapped to each other by the linear automorphism of the
cubic
\begin{gather}z\rightarrow -z,w_2\rightarrow w_2,w_3\rightarrow -w_3.\end{gather} Orbits,
corresponding to different $\mu^2$, are clearly different. Hence
 the family of orbits is parametrized by the non-negative
parameter $\mu^2$.

The Levi forms of the orbits are as follows:

\begin{gather}
\nonumber - \frac{3}{2} \mu |z|^2 + i z \bar{w_2}-i w_2 \bar{z} +
\frac{3}{16} \mu |w_2|^2.
\end{gather}

Then for each $\mu$
the orbits are homogeneous hypersurfaces with
non-degenerate indefinite Levi form.

\it CASE 2. \rm In this case, as was shown in \cite{orbits}, for a
point $(a,b,c), \im b<|a|^2$ we get the following orbits:

\begin{gather*}
\Im w_3 = - 2 \Re z^2 \bar{z} + 4 \Re z \Im w_2 + |\nu| (|z|^2-\Im w_2)^{\frac{3}{2}},\,\im w_2>|z|^2,\nu\in\RR{}.\end{gather*}

Any such orbit is an open smooth part of a real-analytic set

\begin{gather}
\nonumber (\Im w_3 - 4 \Re z \Im w_2 + 2 |z|^2 \Re z)^2 =  \nu^2
(|z|^2-\Im w_2)^{3},
\end{gather}
lying over $P<0$.

Any orbit except the one with $\nu=0$ has two connected
components, corresponding to two $\nu$ with opposite signs. They
can be mapped to each other by the linear automorphism (4) of the
cubic.  Orbits, corresponding  to different $\nu^2$, are clearly
different.  Hence the family of orbits is parametrized by the
non-negative parameter $\nu^2$.

The Levi form in this case equals

\begin{gather}
\nonumber \frac{3}{2} \nu |z|^2 + i z \bar{w_2}-i w_2 \bar{z} +
\frac{3}{16} \nu |w_2|^2.
\end{gather}

 The determinant of the Levi form is $(\frac{9}{32} \nu^2 -
1)$. Hence for $\nu^2 > \frac{32}{9}$ the hypersurfaces are
strictly pseudoconvex; for $\nu^2 = \frac{32}{9}$ the orbit is
Levi-degenerate, the Levi form has one non-zero eigenvalue; for $
\nu^2 < \frac{32}{9}$  the  orbits have indefinite Levi form.

\it CASE 3. \rm  In that case straightforward calculations show
that the values of the vector fields, which form the basis of the
algebra $\mathfrak{g}$, have rank 4 at each point on the cubic and
rank 5 at each point outside the cubic. Hence the cubic is
the only singular orbit of dimension 4. As it was shown in
\cite{orbits}, there are two orbits in that case: $$\im w_2 = |z|^2,
\quad \im w_3 =2\re (z^2 \bar{z})$$ - the cubic, and $$\im w_2 =
|z|^2, \quad \im w_3 \neq 2\re (z^2 \bar{z})$$ - the complement to
the cubic on the cylindric surface $\im w_2=|z|^2$. The second
orbit has two connected components, which can be mapped to each
other by the linear automorphism (4) of the cubic.

To characterize globally the foliation of the space $\CC{3}$,
given by the group $G$, note that the polynomial
$$Q=\Im w_3 - 4 \Re z \Im w_2 + 2 |z|^2 \Re z$$ is also a relative
invariant of  the action (2)-(3) of weight 3. In terms of the relative invariant
polynomials, the  orbits "over the ball" are given by the condition
$$Q^2=\mu^2 P^{3}, \; P>0,$$ the orbits "over the
complement to the ball" are given by the condition $$Q^2=\nu^2
(-P)^{3}, \; P<0,$$ the  orbits "over the sphere" are
given by the condition $$P=Q=0$$ - the cubic, and
$$Q \neq 0, \; P=0$$  - the complement to the cubic. The obtained
description of the foliation is illustrated by the figure.

\begin{figure}
[h]\centering\includegraphics[0,0][10.6 cm, 7.0 cm]{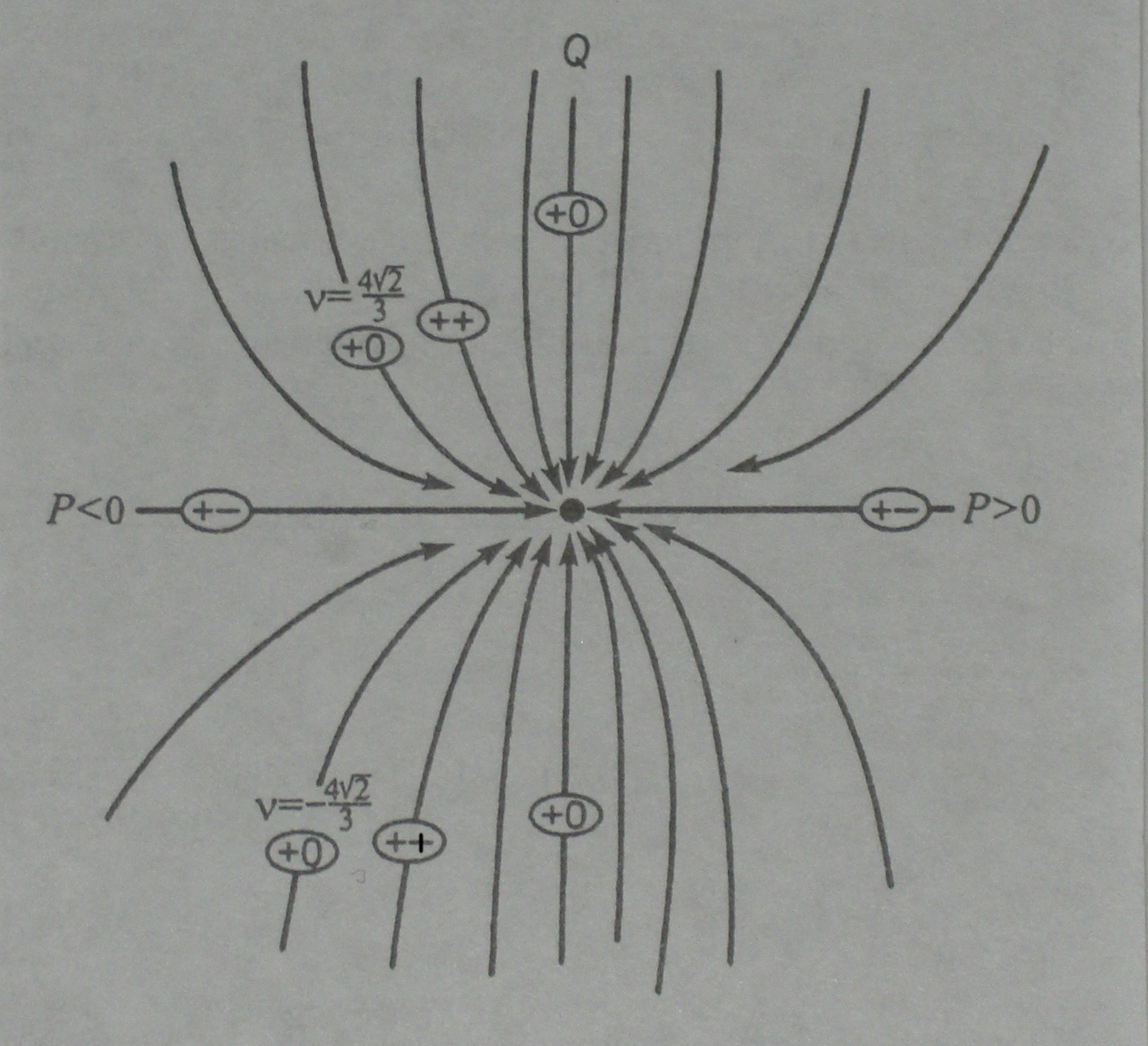}%
\centering\caption {} \label{r:1}%
\end{figure}

Also note the following  fact:  the cubic $C$ is the boundary of any orbit,
and, roughly saying, any two orbits $M,M'$ "meet at $C$", but the
union $M\cup M'\cup Q$ does not form a smooth hypersurface
(moreover, this union does not also decompose to smooth
hypersurfaces), except the case $\mu=\nu=0$, when this union forms the smooth hypersurface
\begin{gather} Q=\Im w_3 - 4 \Re z \Im w_2 + 2 |z|^2 \Re z=0.\end{gather}
Now we  give a tube realization for the obtained foliation in
$\CC{3}$, generated by the group $G$. To do so, remember, that one
of the obtained orbits - corresponding to
$\nu=\frac{4\sqrt{2}}{3}$ - is Levi degenerate with Levi form of
rank 1. For a hypersurface in $\CC{3}$ with Levi form of rank 1 we
have the following dichotomy: it may be either holomorphically
degenerate (in this case it is locally biholomorphally equivalent
to the direct product of a hypersurface in $\CC{2}$ and a the
complex plane, like the 5-dimensional orbit from case 3), or it is
holomorphically non-degenerate and in this case it is
2-nondegenerate (see \cite{bouendy}). It can be checked that our
orbit (denote it by $M$) is 2-nondegenerate. The list of
2-nondegenerate homogeneous surfaces, obtained in \cite{kaup},
consists of one surface with 5-dimensional stabilizer (the tube
over the future light cone), and some surfaces with trivial
stabilizer. Hence $M$ is either isomorphic to the tube over the future light
cone or it has trivial stabilizer and hence is isomorphic to one of
the remaining surfaces in the mentioned list. It is shown in the
next section that $M$ actually has trivial stabilizer and its
infinitesimal automorphism algebra coincides with $\mathfrak{g}$,
so the second possibility holds. It follows from \cite{kaup} that
only one surface in the list - namely the one from example 8.5 - has an infinitesimal
automorphism algebra, isomorphic to $\mathfrak{g}$, which proves,
that $M$ is locally biholomorphically equivalent to the surface
from example 8.5 (denote it by $\tilde{M}$). This surface is a tube over the following affinely
homogeneous hypersurface in $\RR{3}$: $$F=\{c(t)+r
c'(t)\in\RR{3}:r>0,t\in\RR{}\},\,c(t)=(1,t,t^{2}).$$

The infinitesimal automorphism algebra $\widetilde{\mathfrak{g}}$
of $\tilde{M}$ has the following:

\begin{gather*}
\xxx=\dzzz,\,\xx=\dzz,\,\xp=\dz,\\
\x=i\dz+2\z\dzz+3\zz\dzzz\\
\xo=\z\dz+2\zz\dzz+3\zzz\dzzz. \end{gather*}

Since $\tilde{M}$ and $M$ are locally biholomorphically equivalent, there
exists a biholomorphic transformation, defined in a neighborhood
of a point from $\tilde{M}$, which maps this algebra onto $\mathfrak{g}$.
Straightforward calculations show that the mapping

\begin{gather}z=\alpha z_1,w_2=\gamma z_2+\beta z_1^2,w_3=\delta
z_3+\varepsilon z_1^3\end{gather}
with $\alpha=-\frac{i}{\sqrt{2}},\gamma=1,\beta=\frac{i}{2},\delta=\frac{2\sqrt{2}}{3},\varepsilon=\frac{\sqrt{2}}{6}$
indeed maps $\mathfrak{\tilde{g}}$ onto
$\mathfrak{g}$ and hence $\tilde{M}$ onto $M$. This fact gives another possibility to prove that
$M$ is 2-nondegenerate (using the fact that $\tilde{M}$ is
2-nondegenerate). Note that the mapping (6) is a biholomorphic
mapping of $\CC{3}$ onto itself. In particular, it is a global
isomorphism of $\tilde{M}$ and $M$ and the inverse mapping translates all
the orbits from cases 1-3 to some tube homogeneous manifolds in
$\CC{3}$. The corresponding foliation of $\CC{3}$ consists of the
hypersurfaces $N^{\pm}_{\mu},N^0$ (see the introduction) and one
4-dimensional orbit $\tilde{C}=\{y_2=y_1^2,\,y_3=y_1^3\}$. All $N^+_{\mu}$ are Levi-indefinite,
$N^-_{\mu}$ are Levi-indefinite for $\mu<4$, strictly pseudoconvex for $\mu>4$ and 2-nondegenerate for $\mu=4$. The
surface $\tilde{M}$ coincides with $N^{-}_4$ and, unlike all
 other orbits $N_{\mu}^{\pm}$, which are given by equations of degree 6, this orbit is given by
an equation of degree 4 (and of weight 6):

$$y_3^2-3y_1^2y_2^2-6y_1y_2y_3+4y_1^3y_3-4y_2^3=0.$$

  \remark{2.1} It is a very remarkable fact, that the
 mapping (6) transforms the cubic $C$ to the tube $\widetilde{C}$ over the standard twisted
 cubic $$y_2=y_1^2,\,y_3=y_1^3$$ from $\RR{3}$.

\remark{2.2} The approach to the construction of homogeneous
manifolds, used in \cite{orbits}, can be generalized to other
dimensions and model algebras (see \cite{obzor}, \cite{univers}
for the details of the general  notion  of a model manifold) and
can be used as a "machinery" for the construction of homogeneous
CR-manifolds with a "good" Lie transformation group, acting on
them transitively.


\section{Automorphism groups of the orbits and their holomorphic classification }

In this section we  classify the homogeneous hypersurfaces, obtained in the
previous section, up to local biholomorphic equivalence and
compute their automorphism groups. In particular, an analogue of
the Poincare-Alexander theorem is proved for the orbits.

We  parametrize the orbits from cases 1 and 2 by a
non-negative parameter $\mu$ and denote them by $\MP{\mu}$ and
$\MM{\mu}$ correspondingly. Also we  denote the hypersurface
type orbit from case 3 by $M^0$.

To classify the orbits we firstly prove two lemmas.

\lemma{3.1} The infinitesimal automorphism algebra of any orbit
from cases 1,2 is a finite-dimensional algebra of polynomial
vector fields.

\doc All hypersurfaces from cases 1,2 are Levi non-degenerate,
except $\MM{32/9}$.  As it follows from  section 1, the
hypersurface $\MM{32/9}$ is 2-nondegenerate. Hence, according to
\cite{bouendy}, any $M^{\pm}_{\mu}$ has finite-dimensional
infinitesimal automorphism algebra. This algebra contains the
algebra $\mathfrak{g}$ of infinitesimal automorphisms of the
cubic. For each orbit make a translation, which sends the point
$(0,\pm i,i\sqrt{\mu})$ on the orbit to the origin. We  obtain
a surface,  whose infinitesimal automorphism algebra is
finite-dimensional and contains the vector fields $$\dww,\dwww$$
(they come from the translations from $\mathfrak{g}$) and the
vector field
$$z\frac{\partial}{\partial z}+2w_2\dww+3w_3\dwww\pm
2i\dww+3i\sqrt{\mu}\dwww$$ (it comes from the dilation field
$X_0\in \mathfrak{g}$). Hence the new surface contains the origin
and it's complexified infinitesimal automorphism algebra contains
the dilation vector field
$$A=z\frac{\partial}{\partial z}+2w_2\dww+3w_3\dwww.$$ Introducing
weights for the variables and the corresponding weights for the
basic differential operators as in section 2, for a vector field
$X_k$ of weight $k$ we  have $$[A,X_k]=kX_k.$$ Then, expanding any
vector field $X$ from the complexified algebra to a convergent
series $X_{-3}+X_{-2}+X_{-1}+\dots$ near the origin, we get
$$[A,X]=\sum\limits_{k=-3}^{\infty}kX_k.$$ Hence, considering the
minimal polynomial $p$ of the linear operator $\mbox{ad}_A$ on the
complexified algebra, we  have
$$0=p(\mbox{ad}_A)(X)=\sum\limits_{k=-3}^{\infty}p(k)X_k,$$ but $p(k)=0$
only for finite set of integers, hence we get $X_k=0$ for $k$
bigger than some $k_0$, which means that $X$ is polynomial, so the
complexified algebra of the new surface is polynomial, and we can
state the same for the infinitesimal automorphism algebra of the
original surface, as required (see also the remark after the
corollary 4.3 in \cite{kaup}). \qed\\

\lemma{3.2} Suppose that $F$ is a biholomorphic transformation,
which maps a germ of an orbit $M^\pm_\mu$ to a germ of an orbit
$M^\pm_\nu$. Then $F$ is a birational transformation of the
ambient space $\CC{3}$.

\doc In \cite{cubics} the same statement was proved for a
biholomorphic isomorphism $F$ of two germs of cubics. This proof
uses two facts:

1) The infinitesimal automorphism algebras of both surfaces are
finite-dimensional and polynomial.

2) The infinitesimal automorphism algebras of both surfaces
contain vector fields of kind $X_3,\dots,X_0$.

In our case it follows from lemma 3.1 that we can state the same,
hence we obtain the necessary property for $F$, as required. \qed\\

Now we can prove the main statement of this section.

\theor{3.3}

(1) Two orbits $M^\pm_\mu$ and $M^\pm_\nu$ are locally
biholomorphically equivalent if and only if they coincide, except the case $M^+_0\sim M^-_0$, when
both orbits are locally biholomorphically equivalent to the
indefinite quadric $\im z_3 = |z_1|^2-|z_2|^2$ in $\CC{3}$.

(2) All local automorphisms of an orbit $M^\pm_\mu$ belong to $G$
and hence the local automorphism group of $M^\pm_\mu$ coincides
with the identity component of $G$, except the case $\mu=0$, when
the local automorphism group is the image of the identity
component of the 15-dimensional automorphism group of the
indefinite quadric in $\CC{3}$ (see, for example, \cite{obzor})
under a polynomial transformation.

\doc Consider a biholomorphic transformation $F$,
which maps a germ of an orbit $M^\pm_\mu$ to a germ of an orbit
$M^\pm_\nu$, where $\mu> 0$. By lemma 3.2 $F$ is a birational
transformation of the ambient space $\CC{3}$. Denote by $S$ the
singular set of $F$. Since the orbits are holomorphically
non-degenerate, they can not contain an analytic set of dimension
2, hence $M^\pm_\mu/S$ is connected.  Also, since $F$ is rational
and maps a germ of  $M^\pm_\mu$ to a germ of $M^\pm_\nu$, from the
real-analiticity of the orbits we can conclude, that $F$ maps
$M^\pm_\mu/S$ to an open part of $M^\pm_\nu$. Further note, that
the cubic $C$ is generic, so it can not lie in a proper complex
analytic subset of $\CC{3}$, hence there exist an open part of the
cubic such that F is biholomorphic in a neighborhood of this part
(since $F$ is rational). Such a neighborhood contains an open part
of $M^\pm_\mu/S$, because the cubic is the boundary of
$M^\pm_\mu$. This boundary part (since it is essentially singular
for $M^\pm_\mu$, i.e. $M^\pm_\mu$ can not be extended  smoothly to
any neighborhood of any point in the cubic), must go to the
essentially singular (in the above sense) boundary part of
$M^\pm_\nu$. Hence, for $\nu>0$ $F$ must map an open piece of the
cubic to an open piece of the cubic, which implies (see
\cite{cubics}) that $F$ is actually an automorphism of the cubic.
This automorphism preserves all orbits, hence our 2 orbits are
locally biholomorphically equivalent if and only if they coincide,
and in the last case the corresponding biholomorphic automorphism
of a germ of an orbit must belong to the automorphism group of the
cubic. For $\nu=0$ we conclude, that such an $F$ does not exist
(since $M^\pm_0$ has no singular boundary part in the above
sense). So in the case $\nu=0$ the orbits $M^\pm_\mu$ and
$M^\pm_\nu$ are locally biholomorphically inequivalent. It means,
that different orbits $M^\pm_\mu$ and $M^\pm_\nu$ are locally
biholomorphically inequivalent except, may be, the case
$\mu=\nu=0$,  and the automorphism group of a germ of any
$M^\pm_\mu$ for $\mu\neq 0$ coincides with the identity component
of the group $G$. To complete the proof we  show that the
hypersurface (5)  is polynomially equivalent to the indefinite
quadric in $\CC{3}$ (this is sufficient since $M^\pm_0$ are open
parts of this hypersurface, and the quadric is homogeneous).

Considering (5), after a polynomial change of variables, which annihilates the
pluriharmonic terms in the quadratic form $4\re z\im w_2$, we
obtain the following surface: $$\im
w_3=iz\overline{w_2}-i\overline{z}w_2-z^2\overline{z}-\overline{z}^{2}z.$$

The expression in the right side can be presented as
$2\re(z\,\overline{(-iw_2-z^2)})$. So the polynomial
transformation
$$z\rightarrow z,\,w_2\rightarrow -iw_2-z^{2}$$ transforms our
surface to the quadric $\im w_3=2\re(z\overline{w_2})$, which is
linearly equivalent to the indefinite quadric $\im
w_3=|z|^{2}-|w_2|^2$, as required. \qed\\

\cor{3.4} All the orbits $M^\pm_\mu $ for $\mu>0$ have the
property, which is analogue to the Poincare-Alexander theorem for
hyperquadrics: any biholomorphic automorphism of a germ of an
orbit extends to a global automorphism.\rm

\cor{3.5} All the orbits $M^\pm_\mu$  for $\mu>0$ are
holomorphically rigid.

\doc The statement of the corollary follows from the theorem and
the fact, that the group $G$ acts effectively on the orbits from
cases 1,2. \qed\\

It is obvious that the same statements hold also for the tube
manifolds $N^{\pm}_{\mu}$: all $N^{\pm}_{\mu}$ are pairwise
locally biholomorphically  inequivalent except the case
$N^{+}_{0}\sim N^{-}_{0}$. For $\mu>0$ their local automorphisms
turn out to be global, and the local automorphism groups coincide
with the identity component of the image $\widetilde{G}$ of the
group $G$ under the transformation (6). This image is a semidirect
product of the normal subgroup, generated by

\begin{gather*}
\notag z_1\lr z_1+a_1,\,z_2\lr z_2+a_2,\,z_3\lr z_3+a_3,\,a_j\in\RR{}\, -\,\mbox{real translations};\\
\notag z_1\rightarrow z_1+it,\,z_2\rightarrow
z_2+2tz_1+it^2,\,z_3\rightarrow z_3+3tz_2+3t^2z_1+it^3,\,t\in\RR{} \,- \end{gather*}
"translations" along the imaginary direction,
and the subgroup of weighted dilations

\begin{gather}z_1\rightarrow\la z_1,\,z_2\rightarrow\la^2
z_2,\,z_3\rightarrow\la^3 z_3,\,\la\in\RR{*} .\end{gather}

 All $N^{\pm}_{\mu}$ for $\mu>0$
are holomorphically rigid. The manifolds $N^{\pm}_{0}$ are locally
polynomially equivalent to the indefinite quadric in $\CC{3}$.
Their local automorphism groups are 15-dimensional and coincide
with the identity component of the image of the automorphism group
of the indefinite quadric under a polynomial transformation.

\remark{3.6} As well as the claim of remark 2.1, it is a very
remarkable fact that the mapping (6) transforms the automorphism
group $G$ of the cubic to the group $\tilde{G}$, thus giving the
model group $G$ an affine realization.

\section{Local transitive actions of  the model algebra $\mathfrak{g}$
on hypersurfaces in $\CC{3}$}

In the paper \cite{orbits} and in sections 2,3 of  the present paper  the natural action of the
model algebra $\mathfrak{g}$ in the complex space $\CC{3}$ was studied and
  two collections of homogeneous holomorphically
non-degenerate hypersurfaces in $\CC{3}$, on which the algebra
$\mathfrak{g}$ acts transitively, were studied and classified. It is natural to ask now if
 \it all \rm possible transitive actions of this algebra
and \it all \rm possible homogeneous hypersurfaces with
transitively acting Lie algebra $\mathfrak{g}$ have been found. More precisely, it is
natural to formulate the following two problems:

1) To classify all possible local transitive actions of the model
algebra $\mathfrak{g}$ on hypersurfaces in $\CC{3}$ up to local
biholomorphic equivalence.

2) To classify up to local biholomorphic equivalence all locally
homogeneous hypersurfaces in $\CC{3}$, admitting a local transitive
action of the model algebra $\mathfrak{g}$ (\it
$\mathfrak{g}$-homogeneous hypersurfaces). \rm

Clearly, obtaining the first
desired classification, we  reduce the second problem to local
holomorphic classification of the orbits of all possible actions.

We  specify that we call two local holomorphic actions of a
finite-dimensional real Lie algebra $\mathfrak{h}$ on
$\CC{N}_{p_1}$ and  $\CC{N}_{p_2}$ equivalent, if there is a local
biholomorphic mapping $F$ of $\CC{N}_{p_1}$ to  $\CC{N}_{p_2}$,
which translates the first action to the second one, i.e. such
that $\varphi_2\circ\tau=F^{*}\circ\varphi_1$, where
$\varphi_1,\varphi_2$ are the homomorphisms of the algebra
$\mathfrak{h}$ to the algebras of germs of holomorphic vector
fields in the points $p_1,p_2$ correspondingly, $F^*$ is the
natural homomorphism of the algebras of germs of holomorphic
vector fields, induced by $F$, $\tau$ is an automorphism of the Lie algebra
$\mathfrak{h}$. In other words it means, that two realizations of
$\mathfrak{h}$ as an algebra of germs of holomorphic vector fields
are translated to each other by some biholomorphic transformation.
Hence the first classification problem is reduced to the following
one:

\it to classify up to local biholomorphic equivalence all
realizations of the Lie algebra $\mathfrak{g}$ as an algebra of
holomorphic vector fields, defined in a neighborhood of a point
$p\in\CC{3}$, such that their values  at the point
 $p$ (and hence at  any point from a neighborhood of $p$) form a real hypersurface in $\CC{3}$ (and hence in a neighborhood of $p$).\rm

So - we take any algebra of the specified above form, defined in a
neighborhood $U$ of a point $p\in\CC{3}$. Take 5 vector fields
$\xxx,\xx,\x,\xp,\xo$, corresponding by the isomorphism of Lie
algebras to the five basic vector fields from $\mathfrak{g}$,
specified in section 2. Then we  have
the following relations: \begin{align*}&[\xxx,\xx]=0&(32\,)\\
&[\xxx,\x]=0 &(31\,)\\
&[\xxx,\xp]=0 &(31')\\
&[\xxx,\xo]=3\xxx &(30\,)\\
&[\xx,\x]=2\xxx &(21\,)\\
&[\xx,\xp]=0 &(21')\\
&[\xx,\xo]=2\xx &(20\,)\\
&[\x,\xp]=4\xx &(11')\\
&[\x,\xo]=\x &(10\,)\\
&[\xp,\xo]=\xp &(1'0).\end{align*}

Now we  construct a suitable coordinate system for the algebra.
To begin with we  rectify $X_3:\,\,\xxx\lr\dzzz$ - this is
possible since the values of our vector fields have rank 5 in $U$. Let the other fields be:

\begin{gather*}\xx=f_1\dz+f_2\dzz+f_3\dzzz,\x=g_1\dz+g_2\dzz+g_3\dzzz,\\ \xp=h_1\dz+h_2\dzz+h_3\dzzz,
\xo=\lambda_1\dz+\lambda_2\dzz+\lambda_3\dzzz.\end{gather*}
Applying now (32),(31),(31'),(30), we  get: $$\frac{\partial
f_j}{\partial \zzz}=0,\frac{\partial g_j}{\partial \zzz}=0,
\frac{\partial h_j}{\partial \zzz}=0, \frac{\partial
\lambda_1}{\partial \zzz}=0,\frac{\partial \lambda_2}{\partial
\zzz}=0,\frac{\partial \lambda_3}{\partial \zzz}=3.$$ After that
we have two possibilities.

\medskip

1. \it The field $f_1(\z,\zz)\dz+f_2(\z,\zz)\dzz$ is non-zero at
$p$ (general case). \rm Then we rectify this field and have
$$\xx=\dzz+f_3(\z,\zz)\dzzz.$$

(21) gives $\frac{\partial g_1}{\partial \zz}=0,\frac{\partial
g_2}{\partial \zz}=0$. (21') gives $\frac{\partial h_1}{\partial
\zz}=0,\frac{\partial h_2}{\partial \zz}=0$, (20) gives
$\frac{\partial \lambda_1}{\partial \zz}=0,\frac{\partial
h_2}{\partial \zz}=2,$ so now we have
\begin{gather*}\x=g_1(\z)\dz+g_2(\z)\dzz+g_3(\z,\zz)\dzzz,\\
\xp=h_1(\z)\dz+h_2(\z)\dzz+h_3(\z,\zz)\dzzz\\
\xo=\lambda_1(z_1)\dz+(2z_2+\lambda_1(\z))\dzz+(3z_3+\lambda_3(\z,\zz))\dzz.\end{gather*}

Further note, that the equality $g_1=h_1=0$ is impossible, because in that
case the values of our 5 vector fields  have rank $< 5$. So
considering, if necessary, a linear combination $\x+a\xp$
instead of $\x$, which  does not change the relations (32) - (1'0), we
may assume that $g_1\neq 0$ at $p$ and rectify the field $g_1\dz$
(the structure of all other fields  does not change after that),
so now $g_1=1$.

After that, considering  (11'), we  have $\frac{dh_1}{
d\z}=0,\frac{dh_2}{d\z}-h_1\frac{dg_2}{d\z}=4\Rightarrow
h_1=s\in\CC{};h_2=sg_2+4z_1+m.$ Considering (10), we  have
$\frac{d\lambda_1}{dz_1}=1,\lambda_1=\z$ (making a translation
along $\z$ if necessary). Also we  have (from (1'0)):

$\lambda_2'+2g_2-\lambda_1g_2'=g_2;s\lambda_2+2h_2-\lambda_1h_2'=h_2\Rightarrow$
subtracting with the factor  $s$, we get $4z_1+m=4\la_1\Rightarrow
m=0$.

After that we  kill $g_2,h_2$. To do so make the variable change
\begin{gather*}z_2\lr \zz-\int g_2d\z \rr\\
\xxx\lr\xxx,\xx\lr\xx,\\ \x\lr\dz+g_3(\z,\zz)\dzzz,\\ \xp\lr
s\dz+4\z\dzz+h_3(\z,\zz)\dzzz\\
\xo\lr\xo.\end{gather*}

Of course, the functional parameters  change, but their
structure is the same. In the same way after the change

\begin{center}$\zzz\lr\zzz-\int f_3(\z,\zz)dz_2$\,\, we
have\end{center}\,

$\xx\lr\dzz,\xxx\lr\xxx,\x\lr\x, \xp\lr\xp,
\xo\lr\xo.$

Thus after all  transformations

\begin{gather*}
\xxx=\dzzz,\xx=\dzz,\\
\x=\dz+g_3(\z,\zz)\dzzz,\\
\xp=s\dz+4\z\dzz+h_3(\z,\zz)\dzzz,\\
\xo=\z\dz+(2z_2+\la_2(\z))\dzz+(3\zzz+\la_3(\z,\zz))\dzzz.\end{gather*}

Now from (21) we get $\frac{\partial
g_3}{\partial\zz}=2;\,(21')\rr\frac{\partial
h_3}{\partial\zz}=0;\,(20)\rr\frac{\partial\la_3}{\partial\zz}=0$.
As a result we have \begin{gather*}
\xxx=\dzzz,\xx=\dzz,\\
\x=\dz+(2\zz+g_3(\z))\dzzz,\\
\xp=s\dz+4\z\dzz+h_3(\z)\dzzz,\\
\xo=\z\dz+(2z_2+\la_2(\z))\dzz+(3\zzz+\la_3(\z))\dzzz.\end{gather*}

So now we have just one variable functions.

Considering (11'), $h_3'-sg_3'-8z_1=0,h_3=sg_3+4z_1^2+n.$

(10) gives $\la_2'=0,\la_2=0$ (after a translation), and also
$\la_3'+6\zz+3g_3-\z
g_3'-4\zz=2\zz+g_3,\la_3'=z_1g_3'-2g_3,\la_3=z_1g_3-3\int g_3d\z.$

Only one functional parameter  $g_1$ remains, we annihilate it by
the variable change  $\zzz\lr\zzz-\int g_3d\z$, which gives
\begin{gather}
\xxx\lr\dzzz,\xx\lr\dzz,\notag\\
\x=\dz+2\zz\dzzz,\notag\\
\xp=s\dz+4\z\dzz+(4z_1^2+n)\dzzz,\\
\xo=\z\dz+2z_2\dzz+3\zzz\dzzz\notag\end{gather} (the last equality
follows from the formula for $\la_3$ obtained above). Applying also
(1'0), we  get $n=0$ (it follows also from the weights
consideration).

Thus we have a one-parameter collection of polynomial algebras.
Clarify under what assumptions they can be mapped to
$\mathfrak{g}$ - it is not a difficult question now, taking the
polynomiality into account.

Provided we have a biholomorphic mapping of one algebra to another
one, we can state, in particular, that the commutants must be
preserved. It means that
$$\xxx\lr a_3\dwww,\xx\lr a_2\dww+b_2\dwww$$
(we put $w_1:=z$), so

$$\frac{\partial\w}{\partial\zz}=0,\frac{\dd\ww}{\dd\zzz}=0,\frac{\dd\w}{\dd\zzz}=0,\frac{\dd\www}{\dd\zzz}=a_3,
\frac{\dd\ww}{\dd\zz}=a_2,\frac{\dd\www}{\dd\zz}=b_2,$$ that is
$$w_1=F(z_1),w_2=a_2\zz+G(z_1),\www=a_3\www+b_2\ww+H(\z).$$

Also we can state that $\x$ must go to a field from the first
commutant, remembering what such  fields from $\mathfrak{g}$ look
like (see section 2), we get
$$F'(\z)=p,F=p\z$$ (without loss of generality we may assume $F(0)=0$).  Furthermore
$$G'(\z)=2i\w\overline{p}+c_1=2i|p|^2z_1+c_1,G=i|p|^2z_1^2+c_1z_1+c_2$$
and in addition

$ H'+2z_2a_3=i\overline{p}w_1^2+2\re p w_2+c_3=
ip^2\overline{p}\z^2+2\re p(a_2\zz+G(\z))+c_3,$

\noindent which implies $$a_3=a_2\re p.$$

The field $\xp$ goes to the first commutant as well, so firstly we
have $sp=p'$ ($p'$ is the new $p$ for the field $\xp$), further
$$sG'+4a_2z_1= 2iw_1\overline{p'}+d_1=2ip\overline{p'}\z.$$ So
remembering the formula for $G'$ we
get$$2isp\overline{p}=2ip\overline{p'}-4a_2,\rr
a_2=\frac{i}{2}(p\overline{ps}-p\overline{p}s)=|p|^2\im s$$ and
finally
$$sH'+4b_2z_1+4a_3z_1^2=iw_1^2+2\re p'w_2+d_2=ip^2\overline{p'}\z^2+2\re p'(a_2z_2+G(\z))+d_2,$$ so
$ \re p'=\re(ps)=0, a_3=a_2\re p=|p|^2\im s\re p.$

In particular, we see that $\im s\neq 0,\re p\neq 0$. To finish
with $\xp$ it just remains to compare the two obtained formulas for
$H'$. Doing so we get
$$s(ip^2\overline{p}+2i|p|^2\re p)=ip^2\overline{ps}-4a_3\rr-2|pp|^2\im s+2is|p|^2\re p=-4a_3.$$
Applying now the equalities $a_3=|p|^2\im s\re p;\re(ps)=0$ we see
that the obtained above equality holds.

After all calculations we can state that for $\im s=0$ the
necessary transformation is impossible. For all other $s$ we can
take $p=\frac{i}{s},c_i=d_i=0,b_2=0$, and choose $a_2,a_3,G,H$
from the obtained above formulas. All we have to do now is to care
about $\xo$. But one can easily check now that it is sent
exactly to a vector field from $\mathfrak{g}_{0}$.

So we have proved that  for $\im s\neq 0$ we have an equivalence
of $\mathfrak{g}$ and the algebra under consideration. For all
other $s$ the algebras are inequivalent.

Now we clarify, when two algebras with different  $s\in\RR{}$ are
equivalent. Firstly change the field $\xp$ to the field
$\frac{1}{4}(\xp-s\x)$. After that the field $\xp$ has the form:
$$\xp=\z\dzz+(\z^2-2s\zz)\dzzz.$$ All other fields are the same.
After that, taking two algebras for different $s\neq 0,s_1=s,s_2=t,$
make a linear change of variables:
$$w_1=z_1;w_2=\frac{s}{t}z_2;w_3=\frac{s}{t}z_3,$$ then $\xxx,\xx$  dilate, $\x,\xo$ are the same, $\xp$ for $s$  go to $\xp$ for $t$.
It means that such two algebras have the same action in $\CC{3}$.

Thus, in the general case we have 3 algebras: $A (s=i),\, A_0 (s=0), \,A_{1}
(s=1)$. Now we  finally simplify the algebras
$A_0$ and $A_{1}$ (we  suppose $A$ to be simplified as $\mathfrak{\tilde{g}}$).

For $A_1$, putting $s=1$ in (8), after a suitable linear change
we come to the following vector field algebra:

\begin{gather*}
\xxx=\dzzz,\xx=\dzz,\x=\dz+2\zz\dzzz,\notag\\
\xp=s\dz+\z\dzz+z_1^2\dzzz,\xo=\z\dz+2z_2\dzz+3\zzz\dzzz.\notag\end{gather*}

Making the polynomial transformation $$z_1\lr z_1,z_2\lr
z_2-\frac{z_1^2}{2},z_3\lr z_3-\frac{z_1^3}{3},$$ we see that
$X_3\lr X_3,\,X_2\lr X_2, X_1'\lr \dz,\,X_1\lr
\dz-z_1\dzz+2z_2\dzzz,\,X_0\lr X_0$. Finally we have (after a
linear change):
\begin{gather}\xxx=\dzzz,\xx=\dzz,\xp=\dz,\notag\\
\x=\dz+z_1\dzz+2z_2\dzzz,\\
\xo=z_1\dz+2\zz\dzz+3\zzz\dzzz\notag\end{gather}

(of course, the vector fields in (9)  have different from (32) - (1'0) commuting relations, but the algebra, that they generate,
is the same).

For $A_0$ after a linear change we have: \begin{gather}
\xxx=\dzzz,\xx=\dzz,\notag\\
\x=\dz+2\zz\dzzz,\\
 \xp=\z\dzz+z_1^2\dzzz,\notag\\
\xo=\z\dz+2z_2\dzz+3\zzz\dzzz.\notag
\end{gather}

It is shown below that 3 obtained vector field algebras are
inequivalent (it just remains to prove that $A_0$ and $A_1$
are inequivalent).

\rm

\medskip

2. \it The vector field
$f_1(\z,\zz)\dz+f_2(\z,\zz)\dzz$ vanishes at $p$ (degenerate
case). \rm  In that case we rectify
$h_1(\z,\zz)\dz+h_2(\z,\zz)\dzz$ (it's non-zero at $p$ because otherwise
the rank  of the values of our 5 vector fields is less than 5).

After that, applying (21'), (1'0),(11'), we get $\frac{\dd f_3}{\dd
z_1}=0;\frac{\dd \la_1}{\dd z_1}=1;\frac{\dd \la_2}{\dd z_1}=0;
 \frac{\dd g_1}{\dd z_1}=0;\frac{\dd g_2}{\dd z_1}=0$.  Also we can  rectify $g_2(z_2)\dzz$ ($g_2|_p\neq 0$ because of the
 rank). As a result we have
 \begin{gather*}
\xx=f_3(\zz)\dzzz, \x=g_1(\z)\dz+\dzz+g_3(\z,\zz)\dzzz,\\
\xo=(\z+\la_1(\zz))\dz+\la_2(\zz)\dzz+(3\zzz+\la_3(\z,\zz))\dzzz.\end{gather*}

Now (21) gives $-f_3'=2;$ (20) gives $3f_3-\la_2f_3'=2f_3$, so

\begin{center}$f_3=-2z_2+m;\la_2=z_2-m/2.$\end{center}

After a translation $m=0$.  So we have

\begin{center}$\xx=-2\zz\dzzz,\xo=(z_1+\la_1(\z))\dz+\zz\dzz+(3z_3+\la_3)\dzzz.$\end{center}

Making the variable change  $w_1=z_1-\int g_1d\zz$, we  have
$\dzz\lr\dww-g_1\dw$, so $\x\lr\dzz+g_3\dzzz$\,\, and, applying
(10), we get $\la_2'=0,\la_2=0$ (after a translation) and as a
result
$$\xo=\z\dz+\zz\dzz+(3\zzz+\la_3)\dzzz.$$ In the same way, to kill
$h_3$ we make the variable change $w_3=z_3- \int
h_3d\z\rr\dz\lr\dw-h_3\dwww$ and we get \begin{gather*}
\xxx=\dzzz,\,\xx=-2\zz\dzzz,\\
\xp=\dz,\,\x=\dzz+g_3\dzzz,\\
\xo= \z\dz+\zz\dzz+(3\zzz+\la_3)\dzzz.\end{gather*}

After that (1'0) gives $\frac{\dd\la_3}{\dd\z}=0;$ 11' gives
$-\frac{\dd g_3}{\dd\z}=-8\zz\rr\la_3=\la_3(\zz),
g_3=8\z\zz+\varphi(\zz)$. (10) gives $\frac{\dd
\la_3}{\dd\zz}+3g_3-\z\frac{\dd g_3}{\dd\z}-\zz\frac{\dd
g_3}{\dd\zz}=
g_3,\la_3'+16\z\zz+2\varphi-8\z\zz-8\z\zz-\zz\varphi'=0,\la_3'=\zz\varphi'-2\varphi,\la_3=\zz\varphi-3\int\varphi
d\zz.$

It means that $\x=\dzz+(8\z\zz+\varphi)\dzzz$, and after the
variable change $w_3=z_3-\int\varphi
d\zz,\dzz\lr\dww-\varphi\dwww$ we get
$\xo\lr\w\dw+\ww\dww+(-\varphi\ww+3\www+3\int\varphi
dz_2+\ww\varphi-3\int\varphi dz_2)\dwww$ and finally (after a
dilation along $z_3$ and a linear transformation in the algebra)
\begin{gather*}
\xxx=\dzzz,\,\xp=\dz,\\
\xx=\zz\dzzz,\,\x=\dzz+\z\zz\dzzz\\
\xo=\z\dz+\zz\dzz+3\zzz\dzzz.\end{gather*}

We  denote this algebra by $B$. So we have proved that there
are four possible types of local transitive actions of the algebra
$\mathfrak{g}$ on hypersurfaces in $\CC{3}: \,\, A_0,A_1,A,B.$
It is shown in the next section that these four types are
actually inequivalent.

\remark{4.1} Note that the three commuting vector fields $\xxx,\xx,\xp$,
as the case $A_0$ shows, may be linearly dependent over $\CC{}$ at
$p$ and it is impossible to rectify them simultaneously in this
case.

\section{Homogeneous hypersurfaces, associated with the model algebra: explicit presentation, automorphism groups and holomorphic
classification}

In this section we  present the orbits of the obtained holomorphic
vector field algebras $A_0,A_1,A,B$ explicitly, classify the
orbits and compute their infinitesimal automorphism algebras (and
hence the local automorphism groups). It also allows us to prove
the non-equivalence of the algebras $A_0,A_1,A,B$.

Now we  study each of the actions $A_0,A_1,A,B$ separately.

\it CASE $A$. \rm As it was proved in the previous section, the
algebra $A$ is equivalent to the algebras $\mathfrak{g}$ and
$\mathfrak{\tilde{g}}$. So a transitive action of each algebra of
the type $A$  on hypersurfaces in $\CC{3}$ is equivalent to the
action of the algebra $\mathfrak{\tilde{g}}$ near a point
$p\in\CC{3}$,  which satisfies $(\im p_2-(\im p_1)^2)^2+(\im
p_3-(\im p_1)^3)^2>0$. The collection of orbits is
$\{N^+_\mu,N^-_\mu,N^0\},\,\mu\geq 0$. The automorphism groups of
the orbits (and hence the corresponding infinitesimal automorphism
algebras) and their classification were specified in section 3.

\it CASE $A_1$. \rm The vector field algebra (9) (we  also denote
it by $A_1$) acts transitively on hypersurfaces in $\CC{3}$ in a
neighborhood of any point $p\in\CC{3}$ such that $\im p_1\neq 0$.
The corresponding  transformation group  is  a semidirect product
of the normal subgroup, generated by the subgroups

\begin{gather*}
\notag z_1\lr z_1+a_1,\,z_2\lr z_2+a_2,\,z_3\lr z_3+a_3,\,a_j\in\RR{}\, -\,\mbox{real translations};\\
 \notag z_1\lr z_1+t,\,z_2\lr z_2+tz_1+t^2/2,\,z_3\lr
z_3+2tz_2+t^2z_1+t^3/3,t\in\RR{},\end{gather*}
and the subgroup of weighted dilations \begin{gather*}z_1\rightarrow\la z_1,\,z_2\rightarrow\la^2
z_2,\,z_3\rightarrow\la^3 z_3,\,\la\in\RR{*} .\end{gather*}

The foliation to orbits is as specified in the the main theorem.
Also note that all $S_\gamma$ with $\gamma>0$ are linearly
equivalent to $S_1$,  all $S_\gamma$ with $\gamma<0$ are linearly equivalent to $S_{-1}$  by means of the linear transformations
$$z_1\lr  z_1,\,z_2\lr\frac{1}{|\gamma|}\,z_2,\,z_3\lr\frac{1}{\sqrt{|\gamma|}}\,z_3,$$
$S_0$ is locally linearly equivalent to the tube
$$S=\left\{y_3^2=y_1^2+y_2^2,\,y_3>0\right\}$$
over the future light cone (see \cite{kaup1} for more information
about $S$ and $S_0$).

It is easy to see that $S_1$ is strictly pseudoconvex, and
$S_{-1}$ has indefinite Levi form in all points. $S_0$ is Levi
degenerate, more precisely, it is 2-nondegenerate. Hence
$S_1,S_{-1}$ and $S_0$ are locally biholomorphically inequivalent.

Now we  compute the infinitesimal automorphism algebras of $S_1$
and $S_{-1}$ (the infinitesimal automorphism algebra of $S_0$ is
well-known, see \cite{kaup1}).

\proposition{5.1} The infinitesimal automorphism algebras of the
orbits $S_1$ and $S_{-1}$ coincide with the algebra $A_1$, so the
homogeneous hypersurfaces $S_1$ and $S_{-1}$ are holomorphically
rigid.

\doc Our arguments are similar to the proof of lemma 2.1.
Firstly note, that both $S_1$ and $S_{-1}$ are Levi
non-degenerate, hence their infinitesimal automorphism algebras
are finite-dimensional. These two algebras contain $A_1$. Now make
a translation, which sends a point on a surface (say, on $S_1$) to
the origin. In the same way as in lemma 2.1 we conclude that the
complexified algebra $\mathfrak{h}$ of the new surface then
contains the vector field

$$A=z_1\frac{\partial}{\partial z_1}+2z_2\dzz+3z_3\dzzz$$ and hence, by
introducing the corresponding weights as in lemma 2.1, we conclude
the complexified infinitesimal automorphism algebra $\mathfrak{h}$
of the new surface and the infinitesimal automorphism algebra
$\mathfrak{t}$ of $S_1$ are polynomial.

Now taking an arbitrary polynomial $q$ and expanding, using the
polynomiality,  a vector field $X\in \mathfrak{t}$ as
$X_{-3}+X_{-2}+\dots +X_{k_0}$, where each polynomial vector field
$X_j$ has weight $j$, we get (since $A\in\mathfrak{t}$):
$$q(\mbox{ad}_{A})(X)=\sum\limits_{k=-3}^{k_0}q(k)X_k\in \mathfrak{t}.$$

Since the polynomial $q$ is arbitrary, we conclude that each
$X_k\in \mathfrak{t}$. It means, that $\mathfrak{t}$  is a
finite-dimensional graded Lie algebra of kind
$$\mathfrak{t}_{-3}+\mathfrak{t}_{-2}+\dots+\mathfrak{t}_{k_0}.$$

Now we  compute the graded components of the algebra
$\mathfrak{t}$. Any element of $\mathfrak{t}$ is a polynomial
vector field
$$f\dz+g\dzz+h\dzzz,$$ where $f(z),g(z),h(z)$ are polynomials, which satisfy
the tangency condition: $$\im h=3y_1^2\,\im
f(z)+\frac{2y_2}{y_1}\,\im g-\frac{y_2^2}{y_1^2}\,\im f,\,z\in
S_1.$$

Any vector field from $\mathfrak{t}_{-3}$ has the form $a\dzzz$.
From the tangency condition we get $a\in\RR{}$, so
$\mathfrak{t}_{-3}$ coincides with the $(-3)$ - component of
$A_1$. Any vector field from $\mathfrak{t}_{-2}$ has the form
$b\dzz+cz_1\dzzz$. From the tangency condition we get
$b\in\RR{},c=0$, so $\mathfrak{t}_{-2}$ coincides with the $(-2)$
- component of $A_1$. Any vector field from $\mathfrak{t_{-1}}$ has the form
$a\dz+bz_1\dzz+(cz_1^2+dz_2)\dzzz$.  The tangency condition looks as
$$\im(cz_1^2+dz_2)=3y_1^2\im a+2\frac{y_2}{y_1}\im (bz_1)-\frac{y_2^2}{y_1^2}\im a,$$
which follows $c=\im a=\im b=\im d=0,d=2b$ and hence $\mathfrak{t}_{-1}$
 coincides with  $(-1)$ - component of $A_1$. In the same way, from the tangency condition
and and relations of kind $[\mathfrak{t}_i,X_j]\subset
\mathfrak{t}_{i+j},$ applied to a vector field $X_j$ from a
current graded component and an obtained before graded component
$\mathfrak{t}_i$, we conclude, that $\mathfrak{t}_0$ coincides with
the  $0$ - component of $A_1$, and also
$\mathfrak{t}_1=\mathfrak{t}_2=\mathfrak{t}_3=0$. Now  we prove by
induction that $\mathfrak{t}_k=0$ for $k\geq 3$. Since the base is
proved, it is remained to make an induction step, so we suppose
that we have $\mathfrak{t}_j=0$ for $1\leq j\leq k,k\geq 3$. Take
a vector field $X\in \mathfrak{t}_{k+1}$. Then we have $[X,\dz]\in
\mathfrak{t}_k$ and hence $[X,\dz]=0$, which follows that the
coefficients of $X$ do not depend on $z_1$. Also we get
$[X,\dzz]\in \mathfrak{t}_{k-1}$, so $[X,\dzz]=0$ and the
coefficients of $X$ do not depend on $z_2$, and finally
$[X,\dzzz]\in \mathfrak{t}_{k-2}$ and hence $[X,\dzzz]=0$, which
follows that the coefficients of $X$ do not depend on $z_3$. Since
all the coefficients in $X$ consist of monomials of positive
degree, we conclude that $X=0$, so $\mathfrak{t}_k=0$ for $k>0$.
It means that all the graded components of $\mathfrak{t}$ coincide
with the graded components of $A_1$, and hence $\mathfrak{t}=A_1$,
as required. The proof for the case of $S_{-1}$ is the same. \qed\\

\remark{5.2} This proof is a modification of the proof of
proposition 4.2 in \cite{kaup}.

\it CASE $A_0$. \rm The vector field algebra (10)  acts
transitively on hypersurfaces in $\CC{3}$ in a neighborhood of
any point $p\in\CC{3}$ such that $\im p_1\neq 0$. The corresponding transformation group is
a semidirect product
of the normal subgroup, generated by the subgroups

\begin{gather*}
\notag z_1\lr z_1,\,z_2\lr z_2+a_2,\,z_3\lr z_3+a_3,\,a_j\in\RR{}\,-\,\mbox{real translations};\\
\notag z_1\lr z_1+t,\,z_2\lr z_2,\,z_3\lr
z_3+2tz_2,\,t\in\RR{};\\
\notag z_1\lr z_1,\,z_2\lr z_2+rz_1,\,z_3\lr
z_3+rz_1^2,r\in\RR{},\end{gather*} and the subgroup of weighted dilations
\begin{gather*}z_1\rightarrow\la z_1,\,z_2\rightarrow\la^2
z_2,\,z_3\rightarrow\la^3 z_3,\,\la\in\RR{*} .\end{gather*}

The foliation to
orbits is as specified in the the main theorem. Now we
classify the orbits $Q_{\beta}$.

\proposition{5.3} All the orbits $Q_\beta$ are locally
polynomially equivalent to the indefinite quadric in $\CC{3}$.

\doc Making a polynomial transformation, which annihilates the
pluriharmonic terms in the right side of the defining equation of
$Q_{\beta}$, for each $\beta$ we (locally) get the following
surface:
$$y_3=\frac{3\beta}{4}\im(z_1^2\overline{z_1})-\frac{1}{4}\im(z_1\overline{z_2}).$$

The right side of the last equality can be presented as
$$-\frac{1}{4}\im\left(z_1\overline{(3\beta z_1^2+z_2)}\right),$$ so
the invertable polynomial transformation $$z_1\lr z_1,\,z_2\lr
-\frac{3\beta}{4}\,z_1^2-\frac{1}{4}\,z_2$$ transforms our surface
to the quadric $y_3=\im(z_1\overline{z_2})$, which is clearly
linearly equivalent to the standard indefinite quadric
$$y_3=|z_1|^2-|z_2|^2$$ in $\CC{3}$. Proposition is proved. \qed

\it CASE $B$. \rm The vector field algebra, corresponding to B,
acts transitively on hypersurface in $\CC{3}$ in a neighborhood of
any point $p\in\CC{3}$ such that $\im p_1\im p_2\neq 0$. The
corresponding local transformation group is generated by the
following transformation groups:

\begin{gather*}
\notag z_1\lr z_1+a_1,\,z_2\lr z_2,\,z_3\lr z_3+a_3,\,a_j\in\RR{}\,-\,\mbox{real translations};\\
z_1\rightarrow\la z_1,\,z_2\rightarrow\la
z_2,\,z_3\rightarrow\la^3 z_3,\,\la\in\RR{*} \,- \,\mbox{weighted dilations};\\
\notag z_1\lr z_1,\,z_2\lr z_2,\,z_3\lr z_3+rz_2,\,r\in\RR{};\\
\notag z_1\lr z_1,\,z_2\lr z_2+t,\,z_3\lr z_3+tz_1z_2+t^2z_1/2,
t\in\RR{}.\end{gather*}

 The foliation
to orbits is as specified in the the main theorem. So in
case $B$ all orbits are locally linearly equivalent to the
real hyperplane $y_3=0$.

Collecting all obtained results, we can prove the main
theorem.

\doc To prove (1) it remains to prove that $A_1\nsim A_0$
and that $B$ is not equivalent to each of $A$ - actions. The first
claim follows from the fact that any orbit of $A_1$ is locally
non-equivalent to any orbit of $A_0$, the same for the second
claim: all  orbits in $B$ are Levi-flat, all  orbits for
$A$-actions are not Levi-flat.

To prove (2) it remains to prove that no manifold from
case (a) is equivalent to one of the manifolds  from case (b) (the
non-equivalence between  manifolds from the same case was proved
above, the non-equivalence for other pairs of manifolds follows
from the description of the infinitesimal automorphism algebras).
Such equivalence is impossible because all manifolds in cases
(a),(b) are holomorpically rigid, which implies that an
equivalence mapping between two manifolds is an equivalence
mapping between vector fields algebras $A_1$ and $A$, which are
inequivalent (see section 4).

This completely proves the theorem. \qed\\

\bigskip
\small{\obeylines
 Valery K.Beloshapka
 Department of Mathematics
 The Moscow  State University
 Leninskie Gori, MGU, Moscow, RUSSIA
 E-mail: vkb@strogino.ru
 --------------------------------------------------------------------
 Ilya G.Kossovskiy
 Department of Mathematics
The Australian National University
 Canberra, ACT 0200 AUSTRALIA
 E-mail: ilya.kossovskiy@anu.edu.au }

\end{document}